\documentclass[smallextended,times]{article}       

\usepackage{graphicx}

\usepackage{mathptmx}      

\usepackage{amssymb}
\usepackage{amsmath}
\usepackage{amscd}
\usepackage{amsthm}
\usepackage[all]{xy}
\usepackage{manfnt}

\usepackage{centernot}
\usepackage{hyperref}

\usepackage{bm}
\usepackage[OT2,T1]{fontenc}

\newtheorem{theo}{Theorem}[section]
\newtheorem{prop}[theo]{Proposition}
\newtheorem{lem}[theo]{Lemma}
\newtheorem{cor}[theo]{Corollary}

\theoremstyle{definition}
\newtheorem{defi}[theo]{Definition}

\theoremstyle{remark}
\newtheorem{rem}[theo]{Remark}
\newtheorem{ex}[theo]{Example}
\newtheorem{con}[theo]{Convention}

\newcommand \Br {{\rm{Br}}}

\newcommand \Pic {{\rm {Pic}}}
\newcommand \val {{\rm {val}}}
\newcommand \Gal {{\rm{Gal}}}

\newcommand \ev {{\rm{ev}}}
\newcommand \Hom {{\rm {Hom}}}
\newcommand \End {{\rm {End}}}

\newcommand \Aut{{\rm Aut}}

\newcommand \rank{{\rm rank}}

\newcommand \ov{\overline}

\newcommand \Z {{\mathbb Z}}
\newcommand \Q {{\mathbb Q}}
\newcommand \F {{\mathbb F}}

\def\O{{\cal O}}

\newcommand\Frob{{\rm Frob}}

\newcommand\Kum{{\rm Kum}}

\newcommand\beq{\begin{equation} \label}

\def \Im {{\rm {Im\,}}}
\def\fil{{\rm fil}}

\DeclareSymbolFont{cyrletters}{OT2}{wncyr}{m}{n}
\DeclareMathSymbol{\Sha}{\mathalpha}{cyrletters}{"58}

\usepackage{tikz-cd}

\newcommand{\bthe}{\begin{theo}}
\newcommand{\ble}{\begin{lem}}
\newcommand{\bpr}{\begin{prop}}
\newcommand{\bco}{\begin{cor}}
\newcommand{\bde}{\begin{defi}}
\newcommand{\ethe}{\end{theo}}
\newcommand{\ele}{\end{lem}}
\newcommand{\epr}{\end{prop}}
\newcommand{\eco}{\end{cor}}
\newcommand{\ede}{\end{defi}}
\newcommand{\brem}{\begin{rem}}
\newcommand{\erem}{\end{rem}}
\newcommand{\bex}{\begin{ex}}
\newcommand{\eex}{\end{ex}}
\newcommand{\bcon}{\begin{con}}
\newcommand{\econ}{\end{con}}
\begin{document}

\title{Odd torsion Brauer elements and arithmetic of diagonal quartic surfaces over number fields}
\date{}



\author{Evis Ieronymou \footnote{Department of Mathematics and Statistics, University of Cyprus, P.O. Box 20537, 1678, Nicosia,
Cyprus. Email: ieronymou.evis@ucy.ac.cy} }




\maketitle

\begin{abstract}

We use recent advances in the local evaluation of Brauer elements to study the role played by {\it odd} torsion elements of the Brauer group in the arithmetic of diagonal quartic surfaces
over {\it arbitrary} number fields. We show that over a local field if the order of the Brauer element is odd and coprime to the residue characteristic then the evaluation map it induces on the local points is constant.
Over number fields we give a sufficient condition on the coefficients of the equation, which is mild and easy to check, under which the odd torsion does not obstruct weak approximation. We also note a systematic way to produce $K3$ surfaces over $\Q_2$ with good reduction and a non-trivial 2-torsion element of the Brauer group with Swan conductor zero.
\end{abstract}

\maketitle


\section{Introduction}

It is conjectured that the Brauer-Manin obstruction to the Hasse principle and weak approximation is the only one for $K3$ surfaces over number fields \cite{Skconj}.
In this paper we focus on certain aspects of the arithmetic of diagonal quartic surfaces; one of the simplest classes of $K3$ surfaces and one which has been intensively
 studied (see e.g. \cite{SD2}, \cite{Br.ev.dq}, \cite{I10}, \cite{ISZ}, \cite{IS}).
In the case of diagonal quartic surfaces the more relevant
elements of the Brauer group seem to be the elements whose order is a power of $2$. Nonetheless odd order elements of the Brauer group, which are necessarily transcendental in this case,  can provide arithmetic
information as well. For example they can obstruct weak approximation; there are even cases where this obstruction cannot be explained by elements killed by a power of 2. The main purpose of this note is to examine in
some detail the role played by {\it odd} torsion elements of the Brauer group in the arithmetic of diagonal quartic surfaces over {\it arbitrary} number fields.

 Over the rational numbers the arithmetic of diagonal quartic surfaces has been studied in \cite{IS}, and the role of odd torsion of the Brauer group is basically known, for the simple reason that the possibilites of odd
 torsion are quite restricted and well understood. However, going over to arbitrary number fields, things become much more complicated as over an algebraic closure of the ground field the Brauer group of a diagonal
 quartic surface is isomorphic to $\left ( \Q/\Z \right )^{2}$, and some new ideas are needed in order to control the various evaluation maps involved in the Brauer-Manin obstruction.


The Galois module structure of $\Br(\ov X)$, where $X$ is a diagonal quartic surface over a number field is by now well understood (see e.g. \cite{GS}, where the authors treat {\it any} diagonal surface).
In this note we combine this knowledge with some recent advancements in the evaluation of Brauer elements over local fields (\cite{BN}, \cite{I}) as well as some ideas from \cite{IS}. We gather results concerning local evaluation maps in  Section
\ref{section: local stuff}; this section includes results for more general varieties.  Some of the results in Section \ref{section: local stuff} are of independent interest; for example Proposition \ref{exswzer} and also
the discussion before Example \ref{exa: arithm.appl} where we use a recent result of Lazda and Skorobogatov in order to give a systematic way to produce $K3$ surfaces over $\Q_2$ with good reduction and a non-trivial
2-torsion element of the Brauer group which has Swan conductor zero. Moreover, we believe that some of our methods will be relevant in the study of other special classes of $K3$ surfaces for example diagonal degree $2$
$K3$ surfaces; that is, surfaces of the form $w^2=ax^6+by^6+cz^6$
in weighted projective space $\mathbb P(3,1,1,1)$ (cf. \cite{GLN}).


Our first main result is the following.

\textbf{Theorem A }

{\it Let $K/\Q_p$ be a finite extension and let $X/K$ be the diagonal quartic surface given by
$$
X : ax^4+by^4+cz^4+dw^4=0 \quad \ a,b,c,d\in K^*
$$
Let $\mathcal A\in \Br(X)[\ell^n]$ for some odd prime $\ell$ and suppose that $p\neq \ell$.
Then
$${\rm ev}_{\mathcal A,K}:X(K)\to \Br(K)$$ is constant.
}

If $p$ is odd, then the result follows from a good reduction argument and is stated in Proposition \ref{dq1} (i). If $p=2$ this somewhat unexpected result is proved in \textsection 3 where it appears as Theorem \ref{NewTh}. The arguments of \textsection 3 are completely independent of the rest of the paper.
In order to prove Theorem \ref{NewTh} we use the symmetries of supersingular elliptic curves in characteristic $2$. The ideas to control evaluation maps in \textsection 3
 are quite different from the ones used in other parts of the paper which use the notion of Kato's Swan conductor; this is mainly due to the fact that it is a complete mystery what a good model of $X$ over $\O_K$ would be in the case of  characteristic
 $2$ residue field, hence rendering Swan conductor arguments inapplicable.

Let $F/\Q$ be a number field and $\alpha=(a,b,c,d)\in (F^*)^4$. We denote by $X_{\alpha}/F$  the diagonal quartic surface given by

$$
X_{\alpha}:ax^4+by^4+cz^4+dw^4=0.
$$

Our next result (Corollary \ref{thmain}) gives a condition on the coefficients of the defining equation which ensures that the odd part of the Brauer group does not obstruct weak approximation.
See the remarks after Theorem B for the notation used and practical aspects on the verification and inclusivity of the condition.

\medskip

\textbf{Theorem B }

{\it
Suppose that $G_{F,\alpha}\subset I_{F,{\alpha}}$. Then
$$
X_{\alpha}(\mathbb A_F)^{\Br(X_{\alpha})_{\rm{odd}}}=  X_{\alpha}(\mathbb A_F )
$$
}
\brem
\begin{itemize}

\item[(i)]
In \textsection $4.1$ we define the set  $G_{F,\alpha}$ and in \textsection $4.2$ we define the set $I_{F,{\alpha}}$. These two sets are subsets of the set ${\rm Pr}_{\text{odd}}$ of odd rational
primes.
The set $G_{F,\alpha}$ is a finite set that contains all the odd prime divisors of $|\Br(X_{\alpha})/\Br_0(X_{\alpha})|$ and is easily computable.
The set $I_{F,{\alpha}}$ is cofinite in ${\rm Pr}_{\text{odd}}$ and easily computable; it allows one to control the evaluation maps at places above odd primes.

  \item[(ii)]
  Theorem B tells us that under some conditions on the coefficients, the odd torsion part of the Brauer group does not obstruct weak approximation.
  We know that some condition on the coefficients is necessary for the result to hold (\cite[Thm. 12.8]{Preu}, \cite[Thm 1.2]{IS}\, \cite[Thm. 2.3]{IScorr}). The condition $G_{F,\alpha}\subset I_{F,{\alpha}} $
  is quite often satisfied.

\end{itemize}

\erem

The next result follows from the explicit knowledge of the Galois module structure of $\Br(\ov{X_{\alpha}})$; 
see Corollary \ref{sum.coef}.

\medskip

\textbf{Theorem C}
{\it

Suppose that  $$\val_v(a)+\val_v(b)+\val_v(c)+\val_v(d)\not\equiv 0 \mod 4 $$

for a place $v\in \Omega_F$ above an odd prime $p$.

Then $\left (\Br(X_{\alpha})/\Br_0(X_{\alpha}) \right )_{{\rm odd}}$ is a $p$-group.
}
\brem\label{triv.obs}
It is clear that we can use the above in order to produce many examples of diagonal quartic surfaces over arbitrary number fields such that $\left (\Br(X_{\alpha})/\Br_0(X_{\alpha}) \right )_{{\rm odd}}$ is trivial.
\erem

We give some illustrative examples of the computational versatility of our results.
The first example is an application of Proposition \ref{bm2.pr1}.

\bex\label{ex.one.p}

Suppose that $\mathcal A\in \Br(X_{\alpha})$ has order $p^k$ . Let $m$ be the largest ramification index of a
place $v\in\Omega_F$ above $p$ and suppose that $p>4m+1$. Then

$$
X_{\alpha}(\mathbb A_F)^{\mathcal A}=  X_{\alpha}(\mathbb A_F )
$$

\eex

In the next example it is fairly easy to show that $I_{F,{\alpha}}={\rm Pr}_{\text{odd}}$ and hence
there is no need to compute $G_{F,\alpha}$ in order to apply Theorem B.

\bex

Let $F/\Q$ be an extension of degree $n$ and denote by $D$ the product of the odd primes that ramify in $F$. Suppose that (i) every prime divisor of $D$ is at least $n+2$, (ii) $a,b,c,d \in \mathcal O_F$ and (iii) ${\rm Norm}_{F/\Q}(abcd)$ is coprime to $15D$.

 Then
$$
 X_{\alpha}(\mathbb A_F)^{\Br(X_{\alpha})_{\textrm{odd}}} =  X_{\alpha}(\mathbb A_F ).
$$

\eex


The next example is a computational application of Proposition \ref{ord.cons}.

\bex\label{expon}

Suppose that $F/\Q$ is a Galois extension of degree 10. Then every odd torsion element of $\Br(X_{\alpha})/\Br_0(X_{\alpha})$ is killed by $3\cdot5^2\cdot7\cdot 19\cdot 41$.

\eex

Our last example is a combination of Theorem C and Example \ref{ex.one.p}.

\bex

Suppose that
$$\val_v(a)+\val_v(b)+\val_v(c)+\val_v(d)\neq 0 \mod 4 $$
for a place $v\in \Omega_F$ above an odd prime $p$. Let $m$ be
the largest ramification index of a place $w\in\Omega_F$ above p.
 If $p>4m+1$ then
$$
 X_{\alpha}(\mathbb A_F)^{\Br(X_{\alpha})_{{\rm odd}}} =  X_{\alpha}(\mathbb A_F ).
$$

\eex

Let us mention some links of the contents of this paper with other problems.
\begin{itemize}
  \item[(i)] In general, describing the structure of $\Br(X)/\Br_0(X)$ as an abstract group is quite difficult, especially in the presence of transcendental elements.
It is a simple consequence of our results that in many cases we can see that $\Br(X)/\Br_0(X)$ is a $2$-group, just by looking at the coefficients of the defining equation (see Remark \ref{triv.obs}). See also
Example \ref{expon} for a related type of result.

  \item[(ii)] When is the index of a diagonal quartic surface equal to one? This should be governed by the Brauer-Manin obstruction to the existence of zero-cycles of degree $1$ and the results of this paper concerning
      local evaluation maps are very relevant to calculating this obstruction. Let us note that there are known examples of diagonal quartic surfaces of index $1$ and a Brauer-Manin obstruction to the Hasse principle
      \cite{BX}, but there are also cases where there is a Brauer-Manin obstruction to the Hasse principle, there is no Brauer-Manin obstruction to the existence of zero-cycles of degree one,  but it is unclear whether
      the index is $1$ \cite[Rem. 5.20]{Br.bad.red}.

  \item[(iii)]
Consider the following statement where $X/F$ is a $K3$ surface over a number field

\begin{equation} \label{eq:deg.cap}
X(\mathbb A_F)\neq \emptyset \Rightarrow X(\mathbb A_F)^{\Br(X)_{\textrm{odd}}}\neq \emptyset.
\end{equation}

See \cite[\textsection 1]{CV} for more general context and motivation.
It is known that (\ref{eq:deg.cap}) does not hold for $K3$ surfaces in general (\cite[Thm. 1.2]{CN}, \cite[Thm. 1.1]{BV}), however the situation for diagonal quartics is not clear.
 To be more precise, we know that (\ref{eq:deg.cap})  holds when $X$ is the Kummer variety attached to a $2$-covering of an abelian variety \cite[Thm. 3.3]{SZ}, or when
$X$ is a diagonal quartic surface which admits a genus one fibration (see Proposition \ref{apo.Nak}, cf. \cite[Cor. 1.9]{Nak}).
The link of (\ref{eq:deg.cap}) to our results is twofold. Whenever Theorem B applies we get the much stonger statement that $X(\mathbb A_F)^{\Br(X)_{\textrm{odd}}}=X(\mathbb A_F)$. Moreover, whenever Theorem B applies and the rank of $\Pic(X)$ is one, we get genuinely new cases of diagonal quartic surfaces where (\ref{eq:deg.cap}) holds.

\end{itemize}

The organisation of the paper is as follows. In \textsection 2 we have some local considerations and we study the evaluation maps using the notion of Swan conductor.
In \textsection 3 we continue our local considerations but now we concentrate on finite extensions of $\Q_2$. This section has nothing to do with Swan conductors and instead we are using ideas from \cite{IS} utilizing the symmetries of supersingular elliptic curves in characteristic $2$. Note however that if we knew that $X$ acquires good reduction after an extension of degree a power of $2$ then the results of this section would be immediate from the theory of Swan conductors.
In \textsection 4.1 we study the possible odd prime divisors of $\Br(X_{\alpha})/\Br_0(X_{\alpha})$ and define the sets $G_{F,\alpha}$. This subsection consists of simple applications of the known Galois structure of
$\Br(\ov{X_{\alpha}})$. In \textsection 4.2 we define the sets $I_{F,\alpha}$ and combine the previous results.

\section{Evaluation maps over local fields }
\label{section: local stuff}

The Swan conductor for classes in the cohomology of henselian discrete valuation fields was defined by Kato in \cite[\textsection 2]{Kato1}.
Let $K$ be a finite extension of $\Q_p$, with ring of integers $\O_K$. If $\mathcal X/\O_K$ is smooth and proper with geometrically integral fibres then the special fibre induces a notion of Swan conductor for elements of
the Brauer group of the function field of $\mathcal X$.
The Swan conductor is intimately related to the evaluation map ${\rm ev}_{\mathcal A}:\mathcal X(\O_K)\to \Br(K)$ induced by $\mathcal A\in \Br(X)$, where $X$ denotes the generic fibre. These relations are studied in
detail in \cite{BN}.
In turn, understanding these evaluation maps is important in the theory of the Brauer-Manin obstruction.

In this section we will use the notation of \cite[\S 2, \S 3]{I}, where we refer for more details on the notions used. For the convenience of the reader we recall the following.

\begin{itemize}
\item $K$ is a finite extension of $\Q_p$, with ring of integers $\O_K$ and residue field $F$. We denote by $\pi$ a uniformizer of $\O_K$, and by $e$ the absolute ramification index of $K$.

\item
$\mathcal X$ is a faithfully flat, regular, finite type scheme over $\O_K$, with geometrically integral fibres.

\item
$X/K$  is the generic fibre of $\mathcal X/\mathcal{O}_K$.

\item
$Y/F$  is the special fibre of $\mathcal X/\mathcal{O}_K$.

\item $\ell$ is a prime number (the case $\ell=p$ is allowed)

\end{itemize}

We suppose that $X(K)\neq \emptyset$ and we identify $\Br(K)$ with $\Br_0(X)$. We remind the reader that by definition  $\Br_0(X):=\Im(\Br(K)\to \Br(X))$.

Recall that if $W/k$ is a variety over a field $k$, $G$ is an algebraic $k$-group, and $Z\to W$ is a $W$-torsor under $G$ then we have an induced map $\theta_Z:W(k)\to H^1(k,G)$ (see \cite[\textsection 2.2, pg. 22]{Sk}).
We will need the following result which is basically \cite[Lemma 5.12]{Br.bad.red} extended to the case $\ell=p$.
For the definition of $\tau'_{\ell^k}$ see \cite[\textsection 3 and Prop. 2]{I} and for $H^1(F(Y))$ see  \cite[\textsection 2]{I} and the references therein.

\bpr\label{from.Br}
Suppose that $\mathcal X/\mathcal O_K$ is smooth and let $\mathcal A\in\fil_0(\Br(X))$. Suppose that $\tau'_{\ell^m}(\mathcal A)\in H^1(F(Y))$ has order $\ell^n$
and identify $\Br(K)[\ell^n]$ with $H^1(F,\Z/\ell^n)$ via the canonical isomorphism from local class field theory.

 Let $Z\to Y$ be a $Y$-torsor under $\Z/\ell^n$, whose class in
$H^1(F(Y),\Z/\ell^n)$  equals $\tau'_{\ell^m}(\mathcal A)$.
 Then
 $$
 \Im(\ev_{\mathcal A,K})=\Im(\theta_Z)
 $$

 In particular let $\sigma \in Z^1(F,\Z/\ell^n)$ be a cocycle. Then the class $[\sigma]\in H^1(F,\Z/\ell^n) $ lies in the image of $\ev_{\mathcal A,K}$ if and only if the torsor $Z^{\sigma}\to Y$ has an $F$-point.
\epr

\begin{proof}
  This follows from \cite[Prop. 3.1]{BN} and the definition of $\theta_Z$.

  \end{proof}

  \bex\label{curves.bound}
  Suppose that $C/\mathbb F_q$ is a smooth projective geometrically irreducible curve of genus $g$ over a finite field. Denote by $H$ the image of the natural map
  $H^1(\mathbb F_q,\Z/p)\to H^1(C,\Z/p)$.
  Let $Z\to C$ be a $C$-torsor under $\Z/p$, and denote its class in
$H^1(C,\Z/p)$ by $\alpha$.
\begin{itemize}
  \item[(i)] If $\alpha \in H$ then $\theta_Z$ is a constant map.
  \item[(ii)] If $\alpha \notin H$ and $\sqrt{q}+\frac{1}{\sqrt q}>2(p(g-1)+1)$ then $\Im(\theta_Z)=H^1(\mathbb F_q,\Z/p)$.
Indeed under these assumptions $Z$ is geometrically irreducible, and of genus $p(g-1)+1$ by Hurwitz's theorem. The Hasse-Weil bound ensures that $Z^{\sigma}$ has an $\mathbb F_q$-rational point for every
 $\sigma \in Z^1(\mathbb F_q,\Z/p)$ .
\end{itemize}

  \eex

We now concentrate on $K3$ surfaces.

  \bpr\label{exswzer}
  Let $X/K$ be a $K3$ surface with good reduction and $\mathcal A \in \Br(X)[p^t]$. Then the Swan conductor of $\mathcal A$ is divisible by $p$ and bounded above by $m+te-1$, where $m$ is the smallest integer that is
  greater than $\frac{e}{p-1}$.
  Moreover the following are equivalent
  \begin{itemize}
    \item[(1)] The Swan conductor of $\mathcal A$ is zero.
    \item[(2)] The evaluation map ${\rm ev}_{\mathcal A,K}:X(K)\to \Br(K)$ is constant.
    \item[(3)]  The evaluation map ${\rm ev}_{\mathcal A,L}:X(L)\to \Br(L)$ is constant for any extension $L/K$.
  \end{itemize}
  In addition, when $t=1$ and (1) does not hold then the image of ${\rm ev}_{\mathcal A,K}$ consists of $p$ elements.
  \epr
  \begin{proof}
    The fact that the Swan conductor is bounded above by $m+te-1$ is established in the proof of \cite[Prop. 17]{I}.
     Let $Y$ denote the special fibre of a proper smooth model.
     Suppose that the Swan conductor of $\mathcal A$ equals $n>0$, and let ${\rm rsw}_{n,\pi}=(\alpha,\beta)\in \Omega^2_{F(Y)} \oplus \Omega^1_{F(Y)}$ (see \cite[Def. 2.5]{BN} and \cite[Def. 5.3]{Kato1} for the definitions).
     Applying \cite[Thm. 7.1]{Kato1} to the local ring on the model of a closed point in the special fibre, we deduce that both $\alpha$ and $\beta$ are regular on $Y$.
     Note that $Y$ is a K3 surface, in particular $H^0(Y,\Omega_Y^1)=0$ and $\Omega_Y^2\cong \O_Y$ (see e.g. \cite[proof of Prop. 8]{I}). Therefore we deduce from \cite[Thm. B (3) and Rem. 1.2]{BN} that $p$ divides $n$ and the evaluation map takes at least $p$ distinct values and hence is non-constant. In particular 2 implies 1, and we have also established the last assertion of the Proposition.
     As 3 implies 2 trivially it remains to show that 1 implies 3. First we note that if the Swan conductor of $\mathcal A$ (with respect to  $\mathcal X$) is zero, then the Swan conductor of $\mathcal A_L$ (with respect
     to $\mathcal X_{\mathcal O_L}$) is also zero by \cite[Prop. 6.3]{Kato1}. Moreover $Y$ is a $K3$ surface  and hence simply connected. Therefore 1 implies 3 by \cite[Prop 2 and Prop. 6]{I}.
  \end{proof}

\brem\label{nreflem}
If the order of $\mathcal A$ is coprime to $p$ then (1), (2) and (3) of the above lemma hold automatically. This follows from \cite[Prop. 2 and Prop. 6]{I} (cf. \cite[Prop. 2.4]{CTSk}), since the special fibre is simply connected.
\erem

Using the above result we can construct $K3$ surfaces $X/\Q_2$, with good reduction and a two torsion element $\mathcal A\in \Br(X)$, such that the Swan conductor of $\mathcal A$ is zero.
Note that, given an explicit element of $\Br(X)$, it is difficult to calculate its Swan conductor due to the abstract nature of the definition of the Swan conductor.
 If one is given equations for the smooth model, then for $p$-torsion elements (where $p$ is the characteristic of the residue field),
 one can try to sidestep this difficulty by reversing the problem in the following way: First use the formulae in \cite[\textsection 4]{Kato1} to  construct elements in the Brauer group of the function field with
 specified Swan conductor, and then show that the constructed elements actually belong to $\Br(X)$, cf. \cite{Pag}.

  We follow a different strategy, where equations for the smooth model are not used. We start with a known construction of Brauer classes on Kummer surfaces of products of elliptic curves from \cite{SZ12}.
  The elliptic curves are chosen so that we can apply \cite[Thm. 2]{LS}. Let $E,E'$ be elliptic curves with rational two-torsion and $X=\Kum(E\times E')$.
In \cite[\textsection 3]{SZ12} there is an explicit construction of quaternion algebras over the function field of $X$
and an easily checked criterion for whether the quaternion algebras lie in $\Br(X)$.
Using this, it is relatively easy to find examples where both $E$ and $E'$ have good ordinary reduction and $\mathcal A\in \Br(X)$.
By \cite[Thm. 2]{LS} this procedure already gives a systematic way to produce $K3$ surfaces over $\Q_2$ with good reduction and 2-torsion elements of their Brauer groups.
To see whether the Swan conductor of the constructed element is zero we will use Proposition \ref{exswzer}.
Note that for a given explicit element it is a finite calculation to check whether the corresponding evaluation map is constant on $K$-points.

\begin{ex}
 \label{exa: arithm.appl}
Consider the following elliptic curves over $\Q_2$
$$E:\ \ y^2=x(x+1)(x+16) ,\quad \quad E':\ \ y^2=u(u+7)(u-9).$$
Both curves have good ordinary reduction, given by equations over $\F_2$
$$\tilde{E}:\ \ y^2 + xy + y = x^3 + x^2 + x,\quad \quad \tilde{E'}:\ \ y^2 + uy = u^3 + 1.$$
By \cite[Lemma 3.6]{SZ12} we have that
$$\mathcal A=\big((x+1)(x+16),(u+7)(u-9))\big)\in\Br(X)$$
where $X=\Kum(E\times E')$ is the minimal desingularization of the affine surface
$$z^2=x(x+1)(x+16)u(u+7)(u-9).$$

We claim that the value of $\mathcal A$ at a point $P\in X(\Q_2)$ is always trivial.
It suffices to show this when $P$ is a smooth point of the affine surface above.
Let $$P=(x_0,u_0,z_0), \ \ f(x)=(x+1)(x+16), \ \ g(u)=(u+7)(u-9).$$
Recall that a unit in $\Z_2$ is a square if it equals $1 \mod 8$.
If ${\rm val}(x_0)\leq-3$ or ${\rm val}(x_0)\geq7$ (resp. ${\rm val}(u_0)\leq -2$ or ${\rm val}(u_0)\geq 3$) then $f(x_0)$ (resp. $g(u_0)$) is a square in $\Q_2$ and we are done.
For specific values of ${\rm val}(x_0)$ and  ${\rm val}(u_0)$ the possible values of the quaternion algebra can be determined by a finite magma calculation using properties of the Hilbert symbol from \cite[Ch. XIV
\textsection 4]{Ser.Loc}. The claim follows from a Magma calculation \cite{BCP}.

By \cite[Thm. 2]{LS}  we see that $X$ has good reduction and we can conclude by Lemma \ref{exswzer} that $\mathcal A$ has Swan conductor zero.

\end{ex}
\brem
Computational evidence suggest that this might be a general phenomenon; that is to say 2-torsion elements of the Brauer group of Kummer surfaces over $\Q_2$ attached to the product of elliptic curves with rational
$2$-torsion and good ordinary reduction, have Swan conductor zero.
\erem

Our next result concerns diagonal quartic surfaces.

\bpr\label{dq1}
Let $X/K$ be the diagonal quartic surface given by

$$
X:ax^4+by^4+cz^4+dw^4=0 \quad \ a,b,c,d\in K^*
$$
where $K/\Q_p$ is a finite extension with absolute ramification index $e$, and $p$ odd.

 Denote by $\val$ the normalised valuation of $K$ and let $L/K$ be a finite extension.
Let $\mathcal A \in \Br(X)$ have order $n$.
Suppose one of the following.
\begin{enumerate}
  \item[(i)]  $\gcd(n,2p)=1$;

  \item[(ii)]
   $n=p^k$ and $4e<p-1$;

   \item[(iii)]   $n=p^k$, $$\val(a)=\val(b)=\val(c)=\val(d)=0 \mod 2$$ and $2e<p-1$;

     \item[(iv)]  $n=p^k$, $$\val(a)=\val(b)=\val(c)=\val(d)=0$$ and $e<p-1$.

\end{enumerate}

Then $\ev_{\mathcal A, L}$ is a constant map.
\epr
\begin{proof}
 Let $L'/L$ be an extension of degree $m$. By functoriality and the definition of the invariant map in local class field theory we have the following
  commutative diagram

\[
\begin{tikzcd}[row sep=large, column sep=large]
\mathcal X(\mathcal O_L)\arrow{r}{\ev_{\mathcal A,L}}\arrow{d}&\Q/\Z \arrow{d}{\cdot [L':L]} \\
       \mathcal X(\mathcal O_{L'}) \arrow{r}{\ev_{\mathcal A,L'}} &\Q/\Z \\
\end{tikzcd}
\]
 By adding a square root or a fourth root of the uniformiser of $K$ if necessary, it is easy to see from the equation of $X$ that there exists a Galois extension $M/K$ of degree $m\in\{1,2,4,8\}$ such that $X_M$ has good
 reduction. In case (iii) we choose $M$ so that $m=2$ and in case (iv) we choose $M$ to be $K$.
 By \cite[Prop. 8]{I} and Remark \ref{nreflem} we have that $\ev_{\mathcal A,M}$ is constant in all cases. Let $L'$ be the composite of $L$ and $M$. It follows from Proposition \ref{exswzer} (and Remark \ref{nreflem}) that
 $\ev_{\mathcal A,L'}$ is constant. Since $[L':L]$ is a power of $2$ and the order of $\mathcal A$ is odd, we deduce the constancy of $\ev_{\mathcal A,L}$ from the commutative diagram above and the constancy of
 $\ev_{\mathcal A,L'}$.
\end{proof}

We can say something in the case where the uniformiser divides only one of the coefficients.

\bpr
Let $X/K$ be the diagonal quartic surface given by

$$
X:ax^4+by^4+cz^4+dw^4=0 \quad \ a,b,c,d\in K^*
$$
where $K/\Q_p$ is a finite extension and $p$ is odd.

 Denote by $\val$ the normalised valuation of $K$ and suppose that
     $$\val(a)=\val(b)=\val(c)=0, \quad \val(d)\in\{1,2,3\}. $$

Let $\tilde{\mathcal X} /\mathcal O_K$ be the scheme given by the same equation and let $\mathcal X/\mathcal O_K$ be  $\tilde{\mathcal X}/\mathcal O_K$ minus the singular point of its closed fibre.

Let $\mathcal A \in \fil_0( \Br(X))$ have order $p$.

Suppose that
 $$
  \sqrt{q}+\frac{1}{\sqrt q}>2(2p+1)
  $$
where $q=|F|$.
Then $|\Im(\ev_{\mathcal A, K})|\in\{1,p\}$.
\epr

\begin{proof}
 Let $C/F$ be the curve in $\mathbb P^2_F$ given by
$$
ax^4+by^4+cz^4=0
$$ and denote by $V\subseteq \mathbb P^3_F$ the projective cone over $C$. Let $P=(0:0:0:1)\in V(F)$ and set $V'=V-P$. Note that $P$ is a closed point of $\tilde{\mathcal X}$,
 $\mathcal X=\tilde{\mathcal X}-P$ and $V'$ is the closed fibre of $\mathcal X/\mathcal O_K$. Moreover we can see from the equations that $X(K)=\mathcal X(\mathcal O_K)$.

Note that a hyperplane section of $V$ not passing through $P$ is isomorphic to $C$.  By our assumptions on $|F|$, we have that $C(F)$ contains at least two points
which implies that for a point of $V'(F)$ we can find a distinct point of $V'(F)$ such that the line joining them does not pass through $P$. We will use these facts in the sequel.

We apply Proposition \ref{from.Br} with our model being $\mathcal X$.
Let $Z\to V'$ be a $V'$-torsor under $\Z/p$, whose class in
$H^1(F(V'),\Z/p)$  equals $\tau'_{p}(\mathcal A)$.
Suppose that $\Im(\theta_Z)$ consists of more than one point. To complete the proof it suffices to show that $|\Im(\theta_Z)|=p$.  By assumption there exist $Q,R\in V'(F)$ where $Z\to V'$ specialises to distinct cocycles.
\begin{itemize}
  \item [(a)] Suppose that the line joining $Q$ and $R$ does not pass through $P$. Then we can find a hyperplane section containing both $Q$ and $R$ and not passing through $P$.
Denote by $i:S\to V'$ the corresponding closed immersion. Since $\theta_{i^*Z}$ is not constant by construction and $S$ is isomorphic to $C$, it follows from Example \ref{curves.bound}
that $|\Im(\theta_{i^*Z})|=p$. As $\Im(\theta_{i^*Z})\subset \Im(\theta_Z)$ and $|\Im(\theta_Z)| \le p$ this completes the proof in this case.
  \item[(b)] Suppose that the line joining $Q$ and $R$ passes through $P$. Choose a point $Q\neq T\in V'(F)$ such that the line joining $Q$ and $T$ does not pass through $P$.
  Then the line joining $R$ and $T$ does not pass through $P$ either. Hence we can choose a hyperplane section containing both $Q$ (resp. $R$) and $T$ and not passing through $P$.
Denote by $i_1:S_1\to V'$ (resp. $i_2:S_2\to V'$) the corresponding closed immersion. By our constructions either $\theta_{i_1^*Z}$ or $\theta_{i_2^*Z}$ is not constant. We conclude with the same argument as in case (a).
  \end{itemize}

 \end{proof}

 \brem

   It is clear how to modify the argument to prove a similar result in the case that $X$ is a surface, $\mathcal X/\mathcal O_K$ is proper and $Y$ is isomorphic to the projective cone over a smooth projective curve of
   genus $g$. In this situation and when the order of $\mathcal A$ is coprime to $p$, there is a more general result - see \cite[Thm. 6.5]{Br.bad.red}.

 \erem

\section{Evaluation maps over extensions of $\Q_2$}\label{section: ell curves Q2}

The aim of this section is to show the following.

\bthe\label{NewTh}
Let $K/\Q_2$ be a finite extension and let $X/K$ be the diagonal quartic surface given by
$$
X:ax^4+by^4+cz^4+dw^4=0 \quad \ a,b,c,d\in K^*
$$
Let $\mathcal A\in \Br(X)[\ell^n]$ for some odd prime $\ell$.
Then  $${\rm ev}_{\mathcal A,K}:X(K)\to \Br(K)$$ is a constant map.
\ethe

The idea of the proof is the following. After an extension of degree a power of $2$, the surface $X$ admits a rational map of degree a power of $2$ to the Kummer surface attached to the self product of an elliptic curve $E$, which has good supersingular reduction.
To $\mathcal A$ we can associate a Galois module homomorphism, $\phi\in \End_{\Gamma_K}(E_{{\ell^n}})$.
 Understanding the map $\ev_{\mathcal A}$ boils down to understanding the image of $E(K)/\ell^n$, seen as a subgroup of $H^1(K,E_{\ell^n})$ under the induced map $\phi_*\in \End(H^1(K,E_{\ell^n}))$.
 In particular if  $E(K)/\ell^n$ is $\phi_*$-invariant it follows that $\ev_{\mathcal A}$ is constant.
 In order to show that $E(K)/\ell^n$ is $\phi_*$-invariant we proceed as follows. Since $E$ has good reduction and $\ell$ is odd
 it suffices to prove the similar statement for the reduction $A$ of $E$. The crucial point is that $A$ admits two non-commuting automorphisms over the ground field whose square equals multiplication by $-1$.
There are two different arguments depending on the residue of $\ell \mod 4$. If $\ell\equiv 3 \mod 4$ the image of $E(K)/\ell^n$ coincides with the image of an inflation map and is hence $\phi_*$-invariant.
If $\ell\equiv 1 \mod 4$ the action of the Frobenius decomposes as a direct sum of two inequivalent `1-dimensional' modules, and the result follows from this.

We will need the following lemma to take care of the case $\ell\equiv 3 \mod 4$.

\ble\label{voith.matr}

Let $\beta\in \Z[i]$ and let $p$ be a prime with $p\equiv 3 \mod 4$. Let $j=\val_p(\beta-1)$ and suppose that $j\ge 1$. Consider the action of $\Z/p^m$ on $\Z[i]/p^m$ where $1\in \Z/p^m$ acts as multiplication by $\beta$.
Then $H^1(\Z/p^m,\Z[i]/p^m)\cong (\Z/p^t)^2$, where $t={\rm min}\{j,m\}$.

\ele
\begin{proof}
Consider the elements
  $$
  N=\sum_{k=0}^{p^m-1}\beta^k, \quad D=\beta-1
    $$
    It is well known  that

\begin{equation} \label{eq:H1}
 H^1(\Z/p^m,\Z[i]/p^m)=\ker(N)/\Im(D)
\end{equation}
 (see e.g. \cite[Ch. VII, \textsection 4]{Ser.Loc}), where on the right hand side $N$ and $D$ denote the corresponding endomorphisms of $\Z[i]/p^m$.
 Since $\val_p(D)=j>0$ it is easy to check that $\val_p(N)=m$.
These imply that
$$
\ker(N)/\Im(D)\cong (\Z[i]/p^m)/p^j(\Z[i]/p^m)
$$

The result follows from this.

\end{proof}

We will need to know certain facts about the arithmetic of $A$, where $A/\mathbb F_2$ is the supersingular elliptic curve given by
$$
A: y^2+y=x^3+x^2+x+1$$

\ble\label{giaA1}
Let $p$ be an odd prime and let $F/\mathbb F_2$ be a finite extension of even degree.

\begin{itemize}
  \item[(a)] The $p$-primary subgroup of $A(F)$ is of the form $\Z/p^k\times \Z/p^k$ for some $k$.

  \item[(b)] Suppose that the $p$-primary subgroup of $A(F)$ is isomorphic to $\Z/p^k\times \Z/p^k$ for some $k\ge 1$ and let $L/F$ be an extension of degree $p$.
Then the $p$-primary subgroup of $A(L)$ is isomorphic to $\Z/p^{k+1}\times \Z/p^{k+1}$ .
\end{itemize}
\ele
\begin{proof}
\begin{itemize}
  \item[(a)]  Let 
  $\sigma$ and $\rho$ be the automorphisms of $A/\mathbb F_4$ given by
  $$
  \sigma : (x,y)\mapsto (x+1,x+y+\omega)
  $$
  $$
  \rho:(x,y)\mapsto(x+\omega,\omega^2x+y)
  $$
  where $\omega\in \mathbb F_4$ with $\omega^2+\omega+1=0$. We have that $\sigma^2=\rho^2=-1$ and $\sigma\rho=-\rho\sigma$. 

  In order to prove part (a) of the lemma it suffices to establish the following claim.

  {\it Claim.} If $P\in A(F)$ has order $p^m$ then either ${}^{\sigma} P$ or ${}^{\rho} P$ is not a multiple of $P$.

  {\it Proof of Claim}. Suppose that ${}^{\sigma} P=\lambda P$, for some $\lambda \in \Z$. As $\sigma$ is an automorphism we must have that $p$ does not divide $\lambda$.
   Let $Q={}^{\rho} P$. It follows easily that ${}^{\sigma} Q=-\lambda Q$. If $Q$ were a multiple of $P$ we would have that  ${}^{\sigma} Q=\lambda Q$ and hence $2\lambda Q=0$,
   which is a contradiction since $p$ does not divide $\lambda$.


  \item[(b)] By the theory of elliptic curves over finite fields we have the formula
  $$
  |A(\mathbb F_{4^n})|=4^n+1-(2\cdot(-2)^n)=((-2)^n-1)^2
  $$
  Therefore
  $$
  \frac{|A(\mathbb F_{4^{pn}})|}{|A(\mathbb F_{4^n})|}=\left(\sum_{i=0}^{p-1}\alpha^i \right)^2
  $$
  where $\alpha=(-2)^n$. If $\alpha \equiv 1 \mod p$ it is easy to show that $\sum_{i=0}^{p-1}\alpha^i \equiv p \mod p^2$.
  The result follows from this and part (a).
\end{itemize}

\end{proof}

The following is the crucial result concerning the image of $A(F)$ under the Kummer map that we will need.

\ble\label{gia.red}
Let $p$ be an odd prime and let $F/\mathbb F_2$ be a finite extension of even degree.
Suppose that the $p$-primary subgroup of $A(F)$ is of the form $\Z/p^k\times \Z/p^k$ for some $k\ge 1$.
Let $F_m$ denote the fixed field of the kernel of the representation $\Gamma_F\to \Aut(A_{p^{k+m}})$, where $\Gamma_F=\Gal(\ov F/F)$.

\begin{itemize}
   \item[(i)] $[F_m:F]=p^m$, $\Gal(F_m/F)\cong\Z/p^m$.
  \item[(ii)]

The image of the natural map $A(F)/p^m\xrightarrow{\partial} H^1(F, A_{p^m})$ lies in the image of the inflation map
$$
H^1( \Gal(F_m/F),A_{p^m})\xrightarrow{{\rm inf}} H^1(F,A_{p^m}).
$$
Denote by $B_m$ the inverse image of $\Im(\partial)$ in $H^1(  \Gal(F_m/F) ,A_{p^m})$. Let $\phi\in \End_{\Gamma_F}(A_{p^m})$, and denote by
$\phi_*$ the induced endomorphism of  $H^1( \Gal(F_m/F),A_{p^m})$.
Then $\phi_*(B_m)\subset B_m$.

\end{itemize}

\ele
\begin{proof}
  \begin{itemize}
    \item[(i)] Part (i) is clear from Lemma \ref{giaA1}.
    \item[(ii)]

     It follows from the definition of $\partial$ (see e.g. \cite[Ch. VIII, \textsection 2]{Sil1}) that $\Im(\partial)$ is contained in $\Im({\rm inf})$.
    Recall that both $\partial$ and ${\rm inf}$ are injective maps, and so $B_m\cong \left(\Z/p^{\min\{m,k\}}\right)^2$.
    If $m\le k$ then
        $$
  H^1(  \Gal(F_m/F),A_{p^m})\cong\Hom( \Z/p^m, (\Z/p^{m})^2 )\cong (\Z/p^{m})^2
  $$
  and hence $B_m= H^1(  \Gal(F_m/F),A_{p^m})$; so we are done in this case.

  We now assume that $m>k$. Recall the following.
   \begin{itemize}
    \item[(1)] We have $B_m\subset H^1(  \Gal(F_m/F),A_{p^m}) $, and $|B_m|=p^{2k}$.

    \item[(2)] Suppose that $|F|=4^n$ and let $\alpha=(-2)^n$.  Then $|A( F)|=(\alpha-1)^2$.

    \item[(3)]   Let $\sigma$ denote the automorphism of $A/\mathbb F_4$ given by
  $$
  \sigma : (x,y)\mapsto (x+1,x+y+\omega)
  $$

  where $\omega\in \mathbb F_4$ with $\omega^2+\omega+1=0$. We have that $\sigma^2=-1$. 
   By an abuse of notation we also denote by $\sigma$ the endomorphism it induces on $A_{p^m}$.

  \end{itemize}

 Let $s$ denote the image of Frobenius in $ \Gal(F_m/F)$. We will need to distinguish two cases.




{\it Case $p\equiv 3 \mod 4$:} Let $P\in A(\ov{F})$ be an element of order $p^m$.
If ${}^{\sigma} P=aP$ for some $a\in \Z$, then $-P=a{}^{\sigma} P=a^2P$ and so $(1+a^2)P=0$. Therefore $-1$ is a quadratic residue $\mod p$, which is a contradiction. Hence $P$ and ${}^{\sigma}P$ generate
$A_{p^m}$. The matrix representing the image of $s$ in $\End(A_{p^m})$ with respect to this basis is of the form

$$
\begin{pmatrix}
  a&-b \\
  b&a
\end{pmatrix}
$$

for some $a,b\in \Z$.

 By \cite[Prop. 4.11]{W} we have that

 $$
 a \equiv \alpha \mod p^m
 $$
Therefore $\val_p(a-1)=\val_p(\alpha-1)=k$.  As $s$ acts trivially on $A_{p^k}$ we must have that $\val_p(b)\ge k$.
Therefore we can deduce by  Lemma \ref{voith.matr} (take $\beta=a+bi$) that $| H^1( \Gal(F_m/F),A_{p^m})|=p^{2k}$. Hence $B_m= H^1(  \Gal(F_m/F),A_{p^m})$ as the two groups have the same order.

  {\it Case $p\equiv 1 \mod 4$:} As $-1$ is a quadratic residue $\mod p$, we can choose integers $a$ and $b$ such that $a^2+b^2=p$.
  We will now use the endomorphism of $A_{p^m}$ induced by $\sigma \in \Aut(A)$, which we still denote by $\sigma$.
   Let $C_1=\ker((a+b \sigma)^m)\subset A_{p^m}$, and $C_2=\ker((a-b \sigma)^m)\subset A_{p^m}$
   Then we have the following decomposition of $\Gal(F_m/F)$-modules
   $$A_{p^m}=C_1\oplus C_2 $$
  where each summand is isomorphic to $\Z/p^m$ as an abelian group.
If $s$ acts as the same scalar on both factors then we conclude with the same argument with the one used in the case $p\equiv 3 \mod 4$ (taking $b=0$ for the representing matrix).
 Otherwise $s$ acts as multiplication by $\lambda_i$ on $C_i$ $i\in\{1,2\}$, with $\lambda_1\neq \lambda_2 \mod p^m$.
 Since $\lambda_1\neq \lambda_2 \mod p^m$ it follows that each $C_i$ is $\phi$-invariant.

  Let $D_1=H^1(\Gal(F_m/F),C_1)$ and $D_2=H^1(\Gal(F_m/F),C_2)$. We have the following decomposition.
  $$H^1(\Gal(F_m/F),A_{p^m})=D_1\oplus D_2 $$
  Since each $C_i$ is $\phi$-invariant, the induced map $\phi_*$ respects the above decomposition.
 Note also that $\partial$ is $\sigma$-equivariant and that $p=(a+b\sigma)(a-b\sigma)\in \End(A)$. It follows that
  $$B_m=(B_m\cap D_1) \oplus (B_m \cap D_2).$$
  It is known from the cohomology of finite cyclic groups that $D_i$ is cyclic, $i\in\{1,2\}$. The result now follows from the observation that any endomorphism of a cyclic group maps a subgroup to the same subgroup.

  \end{itemize}
\end{proof}

\bpr\label{big.lem.gia.2}
Let $K/\Q_2$ be a finite extension with even residue degree and residue field $F$. Let $E/K$ be an elliptic curve with good reduction such that its reduction
is isomorphic to $A\otimes_{\mathbb F_2} F$.
Let $p$ be an odd prime and suppose that the $p$-primary subgroup of $A(F)$ is of the form $\Z/p^k\times \Z/p^k$ for some $k\ge 1$.
Let $K_m$ be the fixed field of the kernel of the representation $\Gamma_K\to \Aut(E_{p^{k+m}})$, where $\Gamma_K=\Gal(\ov K/K)$.
Then
\begin{itemize}
  \item[(i)] $\Gal(K_m/K)\cong\Z/p^m$
\item[(ii)]
The image of the natural map $E(K)/p^m\xrightarrow{\partial} H^1(K, E_{p^m})$ lies in the image of the inflation map
$$
H^1(\Gal(K_m/K),E_{p^m})\xrightarrow{{\rm inf}} H^1(K,E_{p^m})
$$
Denote by $B_m$ the inverse image of $\Im(\partial)$ in $H^1( \Gal(K_m/K),E_{p^m})$. Let $\phi\in \End_{\Gamma_K}(E_{p^m})$, and denote by
$\phi_*$ the induced endomorphism of  $H^1( \Gal(K_m/K),E_{p^m})$.
Then $\phi_*(B_m)\subset B_m$.
\end{itemize}

\epr
\begin{proof}
\begin{itemize}
  \item[(i)] Recall that for any extension $K'/K$ with residue field $F'$ we have that the reduction map $E(K')[r]\to A(F')[r]$ is injective for any odd $r$ \cite[Ch. VII, Prop. 3.1]{Sil1}.
  Therefore Lemma \ref{giaA1} implies that $K_m$ is the unique unramified extension of $K$ of degree $p^m$. The result follows from this.
In the notation of Lemma  \ref{gia.red} which we use in the second commutative diagram below, $\Gal(K_m/K)$ is actually canonically isomorphic to $\Gal(F_m/F)$.
    \item[(ii)]

Consider the following commutative diagrams where the vertical maps are the reduction maps
\[
\begin{tikzcd}[row sep=large, column sep=large]
 E(K)/p^m  \arrow{r}{\partial}\arrow{d} &   H^1(K, E_{p^m}) \arrow{d}  \\
     A(F)/p^m  \arrow{r}{\partial} &   H^1(F, A_{p^m}) \\
\end{tikzcd}
\]

\[
\begin{tikzcd}[row sep=large, column sep=large]
 H^1( \Gal(K_m/K),E_{p^m})   \arrow{r}{{\rm inf}}\arrow{d} &   H^1(K,E_{p^m}) \arrow{d}  \\
    H^1( \Gal(F_m/F),A_{p^m})  \arrow{r}{{\rm inf}} &   H^1(F, A_{p^m}) \\
\end{tikzcd}
\]

    Since $p$ is odd and $E$ has good reduction all the vertical maps are isomorphisms.

    The result follows from Lemma \ref{gia.red}.

\end{itemize}

\end{proof}

%



{\it Proof of Theorem \ref{NewTh}.}

 It is enough to prove the theorem after an extension of degree a power of 2. Hence we can assume that  $\sqrt[4]{-1},\sqrt{2}\in K$, the residue field has even degree over $\mathbb F_2$, $K$ contains a root of $x^8+6x^4-3$, and that $X/K$ is given by $x^4+y^4=z^4+w^4$. It can be seen by a direct calculation (e.g. using Magma) that the elliptic curve $E:y^2=x^3-x$
 has good reduction over $\Q_2[x]/(x^8+6x^4-3)$ and that the reduction of $E$ is given by the supersingular elliptic curve with equation $$y^2+y=x^3+x^2+x+1$$ over the residue field. Hence the same statements are true over $K$ as well.




  We now need to recall some facts from \cite[\textsection 5.1]{IS}, where we refer for details. Denote by $E^c$ the elliptic curve $y^2=x^3-cx$, where $c\in K^*$.
  To a point in $R \in  X(K)$ we can associate a point $(Q,P)\in (E^{\tau^2}\times E^{\tau^2})(K)$ for a specific $\tau\in K^*$ (this $\tau$ depends on $R$).
  To $\mathcal A$ we can associate a homomorphism $\phi \in \Hom(E^{\tau^2}_{\ell^n},E^{\tau^2}_{\ell^n})$. This $\phi$ defines an element $\mathcal B\in \Br(X)$ with the same image as
  $\mathcal A$ in $\Br(X)/\Br_0(X)$. Evaluating $\mathcal B$ at $R$ gives the element
  $$
  \chi(P)\cup \phi_*(\chi(Q))\in \Br(K)[\ell^n].
  $$
  Here we suppress $r$ from the notation and still denote by $\chi$ the canonical map on the twisted curve $\chi:E^r(K)\to H^1(K,E^r_{\ell^n})$. Moreover
  $$
  \cup: H^1(K,E^{\tau^2}_{\ell^n})\times H^1(K,E^{\tau^2}_{\ell^n})\to \Br(K)
  $$
 is the pairing induced by the Weil pairing $E^{\tau^2}_{\ell^n}\times E^{\tau^2}_{\ell^n}\to \mu_{\ell^n}$.

Therefore it suffices to show that  $\chi(P)\cup \phi_*(\chi(Q))=0\in \Br(K)$. By functoriality it is enough to prove the analogous statement over $K(\sqrt{\tau})$, therefore we may assume that $\tau=1$.
If $E(K)$ has no $\ell$-torsion then $E(K)$ is $\ell$-divisible and so $\chi$ is the zero map and we are done. Otherwise, Lemma \ref{giaA1} (a) (taking $p=\ell$) ensures that the assumptions of Proposition \ref{big.lem.gia.2} hold. Hence $\phi_*(\Im(\chi))\subset\Im(\chi)$ and we conclude by the fact that $\Im(\chi)$ is an isotropic subspace of $H^1(K, E_{\ell^n})$, see \cite[Prop. 4.11]{Poo}. \qed


\section{Diagonal quartic surfaces over number fields}
In this section we let $F/\Q$ be a number field and $\alpha=(a,b,c,d)\in (F^*)^4$. We denote by $X_{\alpha}/F$  the diagonal quartic surface given by

$$
X_{\alpha}:ax^4+by^4+cz^4+dw^4=0 \quad \ a,b,c,d\in F^*.
$$

\subsection{Odd primes dividing $|\Br(X_{\alpha})/\Br_0(X_{\alpha})|$}

The knowledge of the Galois module structure of $\Br(\ov {X_{\alpha}})$ allows us to limit the possible odd primes that divide $|\Br(X_{\alpha})/\Br_0(X_{\alpha})|$.
This will lead to the definition of the sets $G_{F,\alpha}$. For the convenience of the reader we recall in Proposition \ref{apoGS} the result from \cite[\textsection 3.2]{GS} that we need and we refer there for further details.

Let $\ell$ be an odd prime and let $\mathfrak p$ be a prime of $\Q(i)$ not dividing $2\ell$. The principal ideal $\mathfrak p\subset \mathcal O=\mathbb Z[i]$ has a unique generator $\pi$ which is a primary prime  that is
to say $\pi \equiv 1 \mod (1+i)^3$. Denote by $\ov{ \pi}$ the complex conjugate of $\pi$. Note that $\pi$ (and $\ov{\pi}$) has a well defined image in $\mathcal O/\ell^n$.

For the notation used in the next result, recall that $\Gamma_F$ denotes the group $\Gal(\ov F/F)$.
Lemma \ref{co.inv.N} is actually implicit in \cite[Prop. 3.2]{GS}, and is clear from the contents of loc. cit. To be more precise, let $Y=X_{(1,1,1,1)}$. The Galois representation on $H^2(Y_{\bar \Q},\Q_\ell(1))$ is unramified away from $2\ell$ since the Fermat surface $Y$ has good reduction at such primes; this is a standard consequence of smooth and proper base change theorems. From this, one can deduce Lemma \ref{co.inv.N}, cf. \cite[Prop. 2.9 and \textsection 3.2]{GS}. We will provide a different proof of Lemma \ref{co.inv.N}, using the results of \cite{IS}.
\ble\label{co.inv.N}
Let $Y=X_{(1,1,1,1)}$ and let $L$ be the invariant subfield of the action of $\Gamma_{\Q(i)}$ on $\Br(\ov Y)[\ell^n]$. Then $L/\Q(i)$ is unramified at any
$v\in \Omega_{\Q(i)}$ above an odd prime $p\neq \ell$.
\ele
\begin{proof}
  Let $E/\Q$ be the elliptic curve given by $y^2=x^3-x$ which has complex multiplication by $\mathcal O$. Then, $\Br(\ov Y)[\ell^n]$ is isomorphic to
  ${\rm End}(E[\ell^n])/(\mathcal O/\ell^n)$ as a $\Gamma_{\Q(i)}$-module, see \cite[\textsection 3, in particular formulas (5) and (6)]{IS}. Let $T$ denote the completion of $\Q(i)$ at $v$.
  The elliptic curve $E$ has good reduction over $T$, and hence $E[\ell^n]$ is
  an unramified $\Gamma_{T}$-module by the N\'eron-Ogg-Shafarevich criterion. This implies that ${\rm End}(E[\ell^n])/(\mathcal O/\ell^n)$ is unramified as a $\Gamma_{T}$-module
  and hence it follows that $L/\Q(i)$ is unramified at $v$.
\end{proof}

\bco\label{sum.coef}

Suppose that there exists a place $v\in \Omega_F$ above an odd prime $p$ such that

     $$\val_v(a)+\val_v(b)+\val_v(c)+\val_v(d)\not\equiv 0 \mod 4 $$
where $\val_v$ denotes the normalised valuation of $F_v$.

Then $\left (\Br(X_{\alpha})/\Br_0(X_{\alpha}) \right )_{{\rm odd}}$ is a $p$-group.
\eco

\begin{proof}
  Let $\ell$ be an odd prime distinct from $p$.
  By   \cite[Prop. 3.4]{IS} it suffices to show that $(\Br(\ov X_{\alpha})[\ell])^{\Gamma_F}=0$.
Let $v'\in \Omega_{F(i)}$ be above $v$. Note that $e(v'/v)=1$. Let $M$ be the invariant subfield of the action of $\Gamma_{F(i)}$ on $\Br(\ov Y)[\ell]$.
    It follows from  Lemma \ref{co.inv.N} that $M/F(i)$ is unramified at $v'$.
  By \cite[Prop. 3.3]{IS} $\Br(\ov X_{\alpha})[\ell]$ is isomorphic as a $\Gamma_{F(i)}$-module to the twist of $\Br(\ov Y)[\ell]$ by a cocycle $\sigma$ whose class in $H^1(F(i),\mu_4)$ is represented by $(abcd)^{-1}$.
  Note that because $F(i,\sqrt[4]{abcd})/F(i)$ is ramified at $v'$, we can find an element $\beta\in \Gamma_{M}$ such that $\sigma(\beta)=-1$.
  This implies that $(\Br(\ov X_{\alpha})[\ell])^{\beta}=0$ which concludes the proof.
\end{proof}

We will use Proposition \ref{apoGS} in order to control the order of odd torsion elements in $\Br(X_{\alpha})/\Br_0(X_{\alpha})$.

\bpr\label{apoGS} Let $Y=X_{(1,1,1,1)}$ and identify $\Br(\ov Y)[\ell^n]$ with $\mathcal O/\ell^n$.
If $\mathfrak p$ is coprime to $2\ell$ then $\Frob_{\mathfrak p}$ acts as multiplication by $\frac{\pi}{\ov{\pi}}$.
\epr
\begin{proof}
  This is \cite[Prop. 3.2 (i)]{GS}
\end{proof}


We will need some more notation. Denote by $D_{\mathfrak p,m;\ell^n}$ the kernel of multiplication by the image of $(\frac{\pi}{\ov{\pi}})^m-1$ in $\mathcal O/\ell^n$.
Let $A_{\mathfrak p,m;\ell}$ be the least positive integer $\tau$ such that the exponent of $D_{\mathfrak p,m;\ell^{\tau}}$ is less than $\ell^{\tau}$, if such an integer exists, and $\infty$ otherwise.
It is not difficult to see that $A_{\mathfrak p,m;\ell}$ is finite if and only if $\left(\frac{\pi}{\ov{\pi}}\right )^m\neq 1$ if and only if $\mathfrak p$ is above a rational prime which equals $1$ modulo $4$, in which case we have
\begin{equation}\label{form.for.A}
A_{\mathfrak p,m;\ell}={\rm max}\left \{ \val_v(\left(\frac{\pi}{\ov{\pi}}\right )^m-1)   \ \ : \ \ \text{$v\in \Omega_{\Q(i)}$ above $\ell$}             \right \}+1.
\end{equation}


Moreover, it follows from Proposition \ref{apoGS} that  if $k/\mathbb Q(i)$ is a finite extension and $\Frob_{\mathfrak p}^m\in \Gamma_k$ then $\Br(\ov Y)^{\Gamma_k}$ does not contain an element of order
$\ell^{A_{\mathfrak p,m;\ell}}$.

We set

$$\phi(\ell,m):={\rm min}\{ A_{\mathfrak p,m;\ell}-1  \ \ : \ \ \text{$\mathfrak p$ does not divide $2\ell$}\}$$

$$S_m:=\{\ell \in {\rm Pr}_{\text{odd}} \ \ : \ \   \phi(\ell,m)>0  \}$$

$$
G_{F,\alpha}:=S_N$$
Here $N$ is the exponent of $\Gal(F_{\alpha}/\Q(i))$, where $F_{\alpha}$ is the Galois closure of $F(i,\sqrt[4]{abcd})/\Q(i)$.
We use the same notation in the next Proposition. Note that $\phi(\ell,m)$ is always finite.

\bpr\label{ord.cons}
 Let $\ell$ be an odd prime and $n=\phi(\ell,N)$.
Then $\Br(X_{\alpha})/\Br_0(X_{\alpha})$ does not contain an element of order $\ell^{n+1}$.

\epr
\begin{proof}
 Choose a prime  $\mathfrak p$ of $\Q(i)$ not dividing $2\ell$ such that $n=A_{\mathfrak p,N;\ell}-1$. By our assumptions $\Frob_{\mathfrak p}^N\in \Gamma_{F_{\alpha}}$ and so $\Br(\ov Y)^{\Gamma_{F_{\alpha}}}$ does not
 contain an element of order $\ell^{n+1}$. By \cite[Prop. 3.3]{IS} the same is true for $\Br(\ov{ X_{\alpha}})^{\Gamma_{F_\alpha}}$ and a fortiori for $\Br(\ov{ X_{\alpha}})^{\Gamma_{F}}$ as well. We conclude by
 \cite[Prop. 3.4]{IS}.

\end{proof}

\bco\label{poss. prim.}
 Let $\ell$ be an odd prime that divides $|\Br(X_{\alpha})/\Br_0(X_{\alpha})|$.
Then $ \ell\in  G_{F,\alpha}$.

\eco
\begin{proof}
Under the assumption $\Br(X_{\alpha})/\Br_0(X_{\alpha})$ contains an element of order $\ell$. We conclude by Proposition \ref{ord.cons}.
\end{proof}

\bex
Suppose that $F/\Q$ is a Galois extension of degree 10. Then every odd torsion element of $\Br(X_{\alpha})/\Br_0(X_{\alpha})$ is killed by $3\cdot5^2\cdot7\cdot 19\cdot 41$.
\eex
\begin{proof}
  We apply Proposition \ref{ord.cons}. In this case we can take $N=40$ and we use formula (\ref{form.for.A}).
  Note that for specific $\mathfrak p$ and $m$, an easy Magma calculation gives us upper bounds for the various $\phi(\ell,m)$ for $\ell$ coprime to $\mathfrak p$.
  In this case, by considering  $\mathfrak p$ above $5$ we see that an element of odd order coprime to $5$ must have order dividing $3\cdot7\cdot 19\cdot 41\cdot 79\cdot 479\cdot2879$.
  By considering  $\mathfrak p$ above $13$ we see that there are no elements of order $479$ or $2879$ and an element of order a power of $5$ has order at most $25$.
   By considering  $\mathfrak p$ above $17$ we see that there are no elements of order $79$.
\end{proof}


\subsection{Consequences of the local considerations}
For the rest of this section we will abuse notation slightly and often write $X$ instead of $X_{\alpha}$.
Our first task is to define the sets $I_{F,\alpha}$. The idea behind the definition is that we impose conditions to ensure that over the various completions
every element of the Brauer group has Swan conductor zero, at least after an extension of degree a power of $2$.
 We introduce the following  notation.
$$
I_F:=\{p \ \ \text{prime }\ \ :\ \  4e_{w/p}<p-1 \ \ \forall w\in \Omega_F\ \ \text{above $p$}\}.
$$

Denote by $I_{F,\alpha}$ the set of rational primes $p$ such that for all
$v\in \Omega_F$ above $p$ one of the following holds:
\begin{itemize}
  \item[(i)]  $4e_{v/p}<p-1$;
  \item[(ii)] $2e_{v/p}<p-1$ and $\val_v(a)\equiv \val_v(b) \equiv \val_v(c) \equiv \val_v(d) \equiv 0 \mod 2$;
  \item[(iii)]   $e_{v/p}<p-1$ and $\val_v(a)=\val_v(b)=\val_v(c)=\val_v(d)=0 $,
\end{itemize}
  where $\val_v$ denotes the normalised valuation of $F_v$.


\bex
Suppose that $[F:\Q]=n$ and the smallest odd prime that ramifies in $F$ is at least $4n+2$.
Then ${\rm Pr}_{\text{odd}}-\{3,5\}\subset I_{F,\alpha}$.

\eex


\bpr\label{bm2.pr1}

Suppose that $\mathcal A\in \Br(X)$ has order $p^k$ and $p\in I_{F,\alpha} $. Then
$$
X(\mathbb A_F)^{\mathcal A}=
  X(\mathbb A_F )
$$

\epr

\begin{proof}
Clearly there exists an extension $L/F$ of degree a power of $2$ such that $X(L)\neq\emptyset$. Trivially
 $$
 X(\mathbb A_L)^{\Br(X_L)_{\textrm{odd}}}             \neq \emptyset.
  $$
 Since $[L:F]$ is a power of $2$, it would therefore suffice to show that  $\ev_{\mathcal A,L_w}$ and $\ev_{\mathcal A,F_v}$ are constant when $w\in \Omega_L$ and $v\in \Omega_F$.
  Let $v\in \Omega_F$ be above a rational prime $\ell$. There exists $w\in \Omega_L$, above $v$ such that $[L_w:F_v]$ is coprime to $p$. By the diagram in the proof of Proposition \ref{dq1} it therefore suffices to show
  constancy of  $\ev_{\mathcal A,L_w}$ for all $w\in \Omega_L$.

 It follows from the definition of $I_{F,\alpha}$ and Proposition \ref{dq1} that $\ev_{\mathcal A,L_w}$ is constant when $w\in \Omega_L$ is above an odd prime $\ell$. If $w$ is above $2$ then $\ev_{\mathcal A,L_w}$
  is constant because of Theorem \ref{NewTh}.
\end{proof}

\bco\label{thmain}

Suppose that $G_{F,\alpha}\subset I_{F,{\alpha}} $. Then
$$
X_{\alpha}(\mathbb A_F)^{\Br(X)_{{\rm odd}}}=  X_{\alpha}(\mathbb A_F )
$$
\eco
\begin{proof}
   Proposition \ref{bm2.pr1} and Corollary \ref{poss. prim.}.
\end{proof}

\ble

Suppose that $X(\mathbb A_F)\neq \emptyset$ and that $\mathcal A_i\in \Br(X)$ has order $p_i$ for $i\in[1,...n]$, where the $p_i$ are pairwise distinct odd primes.
Suppose that for every place $v\in\Omega_F$ above any $p_i$, $X_v$ has good reduction.

Then

$$
X(\mathbb A_F)^{B}\neq \emptyset
$$
where $B$ is the subgroup of $\Br(X)$ generated by the $\mathcal A_i$.

\ele
\begin{proof}
Fix $i\in[1,...,n]$. By Theorem \ref{NewTh} we know that $\ev_{\mathcal A_i,F_v}$ is constant for any $v$ above $2$.
Moreover from Proposition \ref{dq1} it follows that $\ev_{\mathcal A_i,F_v}$ is constant for any $v$ above $\ell$, when $\ell \neq p_i$ (and all extensions of $F_v$).
Suppose that $\ev_{\mathcal A_i,F_v}$ is constant for any $v$ above $p_i$. It follows from Proposition \ref{exswzer} that the same is true for any extension of $F_v$ as well.
Since there exists $L$ such that $[L:F]=4$ and $X(L)\neq\emptyset$ it follows that
  $
 X(\mathbb A_F)^{\mathcal A_i}=      X(\mathbb A_F)
  $. Therefore we may assume that for all $i$ there exists a place $v_i$ above $p_i$ such that $\ev_{\mathcal A_i,F_{v_i}}$ is non-constant.

  Let $z=(z_v)_{v\in \Omega_F}$ be a family of local rational points. From what we showed before if we change $z_{v_i}$ the value of the sum of the evaluation maps at $\mathcal A_j$ is unaffected for $j\neq i$.
  Moreover since by  Proposition $\ref{exswzer}$ we know that $\ev_{A_i,F_{v_i}}$ takes $p_i$ values we can change $z_{v_i}$ in such a way so that the sum of the evaluation maps at $\mathcal A_i$ is zero.
  Since we can do this separately for each $i$ the result follows from the last two observations.
\end{proof}

\subsection{A result for fibrations }

We record the following application of \cite[Thm. 1.4]{Nak}) which is concerned with the case when $X_{\alpha}$ admits a genus one fibration.
\bpr\label{apo.Nak}
Suppose that $X_{\alpha}$ admits a genus one fibration. Then

$$
X_{\alpha}(\mathbb A_F)\neq \emptyset \Rightarrow X_{\alpha}(\mathbb A_F)^{\Br(X_{\alpha})_{{\rm odd}}}\neq \emptyset.
$$

\epr
\begin{proof}
Let $X_{\alpha}\to \mathbb P^1$ be a genus one fibration with generic fibre $X_{\eta}$. Let $l$ denote the index of $X_{\eta}$. By \cite[Thm. 1.4]{Nak}) it suffices to show that $l$ is a power of $2$.
It is enough to prove the last statement after we make an extension of the ground field of degree a power of $2$. Therefore we assume that $\Pic(X_{\alpha})=\Pic(\ov{X_{\alpha}})$.
Let $T$ denote the class of a general fibre in $\Pic(X_{\alpha})$. The multisection index $m$ of the fibration is either 1 or 2; this is because $T/m$ belongs to $\Hom(\Pic(X_{\alpha}),\Z )$ and induces an element of
order $m$ in the discriminant group of $\Pic(X_{\alpha})$, see \cite[\textsection 2 and pg. 2081 Example]{Keum}. We conclude by noting that $l=m$; as follows from \cite[Ch. 9, Prop. 1.30]{Liu}

\end{proof}

\brem\label{cond.ge.1.fibr}
The following are well-known. See e.g. \cite[Cor. 2.27]{Br02} for part (i) and \cite[\textsection 1]{SD2} for part (ii); the proofs given there work for arbitrary number fields.
Suppose that $X_{\alpha}(\mathbb A_F)\neq \emptyset$. Then the following hold.

\begin{itemize}
  \item[(i)]
   $X_{\alpha}$ admits a genus one fibration if and only if one of the following holds
\begin{itemize}
  \item[(1)] $\rank(\Pic(X_{\alpha}))\geq 3$
  \item[(2)]  $\rank(\Pic(X_{\alpha}))=2$ and the determinant of the intersection form on $\Pic(X_{\alpha})$ is minus a square.

\end{itemize}

  \item[(ii)] If $abcd$ is a square in $F$ then $X_{\alpha}$ admits a genus one fibration.
\end{itemize}

\erem

\begin{center}
   ACKNOWLEDGMENTS
\end{center}

 The author is grateful to Alexei Skorobogatov and Alp Bassa for helpful discussions related to the contents of this paper.
Thanks are also due to the anonymous referee for numerous helpful comments.

\bibliographystyle{amsplain}

\begin{thebibliography}{99}

\bibitem{BV} J. Berg and A. Várilly-Alvarado, \textit{Odd order obstructions to the Hasse principle on general K3 surfaces.}
Math. Comp. {\bf 89} (2020), no. 323, 1395-1416.


\bibitem{BCP} W. Bosma, J. Cannon, and C. Playoust,  \textit{The Magma algebra system. I. The user language}, J. Symbolic Comput., {\bf 24} (1997), 235–265.



\bibitem{BX} A. Bremner and T.N. Xuan,  \textit{An interesting quartic surface, everywhere locally solvable, with cubic point but no global point.}
Publ. Math. Debrecen {\bf 93} (2018), no. 1-2, 253-260.


\bibitem{Br02} M. Bright,  \textit{Computations on diagonal quartic surfaces.} Ph.D. dissertation. University of Cambridge, Cambridge, 2002.


\bibitem{Br.ev.dq} M. Bright,  \textit{The Brauer-Manin obstruction on a general diagonal quartic surface.}
Acta Arith. {\bf 147} (2011), no. 3, 291-302.

\bibitem{Br.bad.red} M. Bright,  \textit{Bad reduction of the Brauer-Manin obstruction.}
J. Lond. Math. Soc. (2) {\bf 91} (2015), no. 3, 643-666.


\bibitem{BN}  M. Bright and R. Newton,  \textit{Evaluating the wild Brauer group.} preprint arXiv:2009.03282.




\bibitem{CTSk} J.-L. Colliot-Th\'{e}l\`{e}ne, and A.N. Skorobogatov,  \textit{Good reduction of the Brauer-Manin obstruction.}
Trans. Amer. Math. Soc. {\bf 365} (2013), no. 2, 579-590.





\bibitem{CN} P. Corn and M. Nakahara,  \textit{Brauer-Manin obstructions on degree 2 K3 surfaces.}
Res. Number Theory {\bf 4} (2018), no. 3, Paper No. 33, 16 pp.


\bibitem{CV} B. Creutz and B. Viray,  \textit{Degree and the Brauer-Manin obstruction.}
With an appendix by A.N. Skorobogatov,
Algebra Number Theory {\bf 12} (2018), no. 10, 2445-2470.


\bibitem{GS} D. Gvirtz and A.N. Skorobogatov,  \textit{Cohomology and the Brauer groups of diagonal surfaces.}
Duke Math. J. {\bf 171} (2022), no. 6, 1299-1347.


\bibitem{GLN} D. Gvirtz, D. Loughran and M. Nakahara,  \textit{Quantitative arithmetic of diagonal degree 2 K3 surfaces.} Math. Ann. {\bf 384} (2022), no. 1-2, 135-209.


\bibitem{I10} E. Ieronymou,  \textit{Diagonal quartic surfaces and transcendental elements of the Brauer groups.}
J. Inst. Math. Jussieu {\bf 9} (2010), no. 4, 769-798.

\bibitem{I} E. Ieronymou,  \textit{Evaluation of Brauer elements over local fields.} Math. Ann. {\bf 382} (2022), no. 1-2, 239-254.

\bibitem{IS} E. Ieronymou and A.N. Skorobogatov,  \textit{Odd order Brauer-Manin obstruction on diagonal quartic surfaces.}
Adv. Math. {\bf 270} (2015), 181-205.

\bibitem{IScorr} E. Ieronymou and A.N. Skorobogatov,  \textit{Corrigendum to "Odd order Brauer-Manin obstruction on diagonal quartic surfaces'' [Adv. Math. 270 (2015) 181–205] [MR3286534].}
Adv. Math. {\bf 307} (2017), 1372-1377.

\bibitem{ISZ} E. Ieronymou, A.N. Skorobogatov and Y. Zarhin,  \textit{On the Brauer group of diagonal quartic surfaces.}
With an appendix by Peter Swinnerton-Dyer.
J. Lond. Math. Soc. (2) {\bf 83} (2011), no. 3, 659-672.


\bibitem{Kato1} K. Kato,  \textit{Swan conductors for characters of degree one in the imperfect residue field case.} Algebraic K-theory and algebraic number theory (Honolulu, HI, 1987), 101-131,
Contemp. Math., 83, Amer. Math. Soc., Providence, RI, 1989.

\bibitem{Keum} J. Keum,  \textit{A note on elliptic $K3$ surfaces.} Trans. Amer. Math. Soc. {\bf 352} (2000), no. 5, 2077-2086.


\bibitem{Liu} Q. Liu,  \textit{Algebraic geometry and arithmetic curves.}
Translated from the French by Reinie Erné. Oxford Graduate Texts in Mathematics, 6. Oxford Science Publications. Oxford University Press, Oxford, 2002. xvi+576 pp.

\bibitem{LS} C. Lazda and  A.N. Skorobogatov, \textit{Reduction of Kummer surfaces modulo 2 in the non-supersingular case.} preprint arXiv:2205.13831.


\bibitem{Nak} M. Nakahara,  \textit{Index of fibrations and Brauer classes that never obstruct the Hasse principle.}
Adv. Math. {\bf 348} (2019), 512-522.


\bibitem{Pag} M. Pagano,  \textit{An example of a Brauer-Manin obstruction to weak approximation at a prime with good reduction.}
Res. Number Theory {\bf 8} (2022), no. 3, Paper No. 63, 15 pp.

\bibitem{Poo}
B. Poonen and E. Rains,  \textit{Random maximal isotropic subspaces and Selmer groups.}
J. Amer. Math. Soc. {\bf 25} (2012), no. 1, 245-269.

\bibitem{Preu} T. Preu, \textit{Example of a transcendental 3-torsion Brauer-Manin obstruction on a diagonal quartic surface.} Torsors, étale homotopy and applications to rational points, 447-459,
London Math. Soc. Lecture Note Ser., 405, Cambridge Univ. Press, Cambridge, 2013.

\bibitem{Ser.Loc} J.-P. Serre, \textit{Local fields.}
Translated from the French by Marvin Jay Greenberg. Graduate Texts in Mathematics, 67. Springer-Verlag, New York-Berlin, 1979. viii+241 pp.


\bibitem{Sil1} J. Silverman, \textit{The arithmetic of elliptic curves.}
Second edition. Graduate Texts in Mathematics, 106. Springer, Dordrecht, 2009. xx+513 pp.





\bibitem{Sk} A.N. Skorobogatov. \textit{Torsors and rational points.}
Cambridge Tracts in Mathematics, 144. Cambridge University Press, Cambridge, 2001. viii+187 pp.

\bibitem{Skconj} A.N. Skorobogatov, \textit{Diagonal quartic surfaces.} Oberwolfach Rep. {\bf 33} (2009): 76-9.



\bibitem{SZ12} A.N. Skorobogatov and Y. Zarhin, \textit{The Brauer group of Kummer surfaces and torsion of elliptic curves.} J. reine angew. Math. {\bf 666} (2012) 115--140.


\bibitem{SZ} A.N. Skorobogatov and Y. Zarhin, \textit{Kummer varieties and their Brauer groups.}
Pure Appl. Math. Q. {\bf 13} (2017), no. 2, 337-368.



\bibitem{SD2} P. Swinnerton-Dyer, \textit{Arithmetic of diagonal quartic surfaces. II.}
Proc. London Math. Soc. (3) {\bf 80} (2000), no. 3, 513-544.


\bibitem{W} L.C. Washington, \textit{Elliptic curves.}
Number theory and cryptography. Second edition. Discrete Mathematics and its Applications (Boca Raton).

\end{thebibliography}

\end{document}